%&amstex
\input amstex
\documentstyle{amsppt}
\magnification=\magstephalf \hsize = 6.5 truein \vsize = 9 truein
\vskip 3.5 in
%\parindent = 0 pt
%\pageno=2
%\parskip=24 pt

%\NoRunningHeads
\NoBlackBoxes
\TagsAsMath

\def\label#1{\par%
        \hangafter 1%
        \hangindent .75 in%
        \noindent%
        \hbox to .75 in{#1\hfill}%
        \ignorespaces%
        }

%Section macros
\newskip\sectionskipamount
\sectionskipamount = 24pt plus 8pt minus 8pt
\def\sectionskip{\vskip\sectionskipamount}
\define\sectionbreak{%
        \par  \ifdim\lastskip<\sectionskipamount
        \removelastskip  \penalty-2000  \sectionskip  \fi}
\define\section#1{%
        \sectionbreak   %Encourages a page break, else inserts 24pt glue
        \subheading{#1}%
        \bigskip
        }

%QED box, from the TeXbook, p. 106.
\redefine\qed{{\unskip\nobreak\hfil\penalty50\hskip2em\vadjust{}\nobreak\hfil
    $\square$\parfillskip=0pt\finalhyphendemerits=0\par}}

%Other math symbols
        
        \let    \< = \langle
        \let    \> = \rangle

%Operator name macros
\define\op#1{\operatorname{\fam=0\tenrm{#1}}} %for text in math mode use
                                              %\op{...}
 %Greek letters
        \define         \a              {\alpha}
        \redefine       \b              {\beta}
        \redefine       \d              {\delta}
        \redefine       \D              {\Delta}
        \define         \e              {\varepsilon}
        \define         \E              {\op {E}}
        \define         \g              {\gamma}
        \define         \G              {\Gamma}
        \redefine       \l              {\lambda}
        \redefine       \L              {\Lambda}
        \define         \n              {\nabla}
        \redefine       \var            {\varphi}
        \define         \s              {\sigma}
        \redefine       \Sig            {\Sigma}
        \redefine       \t              {\tau}
        \define         \th             {\theta}
        \redefine       \O              {\Omega}
        \redefine       \o              {\omega}
        \define         \z              {\zeta}
        \define         \k              {\kappa}
        \redefine       \i              {\infty}
        \define         \p              {\partial}
        \define         \vsfg           {\midspace{0.1 truein}}

\topmatter
\title Primitivity preserving endomorphisms of free groups
\endtitle
\author  Donghi Lee
\endauthor

\address {Department of Mathematics, University of Illinois at
Urbana--Champaign, 1409 West Green Street, Urbana, IL 61801, USA}
\endaddress

\email {d-lee9\@math.uiuc.edu}
\endemail

\subjclass Primary 20F05, 20F06, 20F32
\endsubjclass

\abstract{An endomorphism $\phi$ of a free group is called
primitivity preserving if $\phi$ takes every primitive element to
another primitive. In this paper we prove that every
primitivity preserving endomorphism of a free group of a finite
rank $n \ge 3$ is an automorphism.}
\endabstract
\endtopmatter

%09/20/00

\document
\baselineskip=24pt

\heading 1. Introduction
\endheading

Let $\phi$ be an endomorphism of a free group. Following [2] we
call $\phi$ a {\it primitivity preserving endomorphism} if $\phi$
takes every primitive element to another primitive.

Let $F_n$ be the free group of a finite rank $n$ on the set ${\frak X}=\{x_1, \dots, x_n\}$. In [1, 5], Shpilrain asked whether every
primitivity preserving endomorphism $\phi$ of $F_n$ is actually an
automorphism, and this problem was solved in the affirmative for
$n=2$ by Ivanov [2] and Shpilrain [6]. Ivanov [2] also proved that
the answer is positive for an arbitrary rank $n \ge 3$ under the
additional assumption that $\phi (F_n)$ contains a
primitive pair (i.e., a pair $\psi(x_1)$, $\psi(x_2)$ for some
$\psi \in \text {Aut} \, F_n$). The purpose of this paper is to
present a positive solution to this problem for rank $n \ge 3$
without Ivanov's extra assumption:

\proclaim {Theorem} Let $F_n$ be the free
group of rank $n \ge 3$ on the set ${\frak X}=\{x_1, \dots, x_n\}$ and $\phi$ a primitivity preserving
endomorphism of $F_n$. Then $\phi$ is an automorphism of $F_n$.
\endproclaim

The proof of the Theorem is developed in a series of several steps.
In Steps A--C, the following result of McCool [4] (also see [3])
plays a crucial role:

\proclaim {Lemma {\rm [4]}} Let $F_n$ be a free group of a finite rank $n$
and $\psi \in \text {Aut} \, F_n$, and let $U \in F_n$ be such
that $\|\psi (U)\| < \|U\|$. Then there exist Whitehead
automorphisms $\a_1, \dots, \a_r$ of $F_n$ such that $\psi = \a_r
\cdots \a_1$ with $\|\a_i \cdots \a_1(U)\| < \|U\|$ for $1 \le i
\le r$.
\endproclaim

Here, by $\|X\|$ we mean the length of a reduced word equal in
$F_n$ to $X$ which is denoted by $\overline X$, whereas $|X|$
denotes the length of $X$ (e.g., $\|x_1x_1^{-1}\|=0$ and
$|x_1x_1^{-1}|=2$). Recall that a Whitehead automorphism $\o$ of
$F_n$ is an automorphism of one of the following two types (see
[2, 7]):

\roster
\item"(W1)" $\o$ permutes elements in ${\frak X}^{\pm 1}$.
\item"(W2)" $\o$ is defined by a set ${\Cal S} \subset {\frak
X}^{\pm 1}$ and a letter $a \in {\frak X}^{\pm 1}$ with $a \in
{\Cal S}$ and $a^{-1} \notin {\Cal S}$ in such a way that if $x
\in {\frak X}^{\pm 1}$ then (a) $\o(x)=x$ provided $x=a^{\pm 1}$;
(b) $\o(x)=xa$ provided $x \neq a$, $x \in {\Cal S}$ and $x^{-1}
\notin {\Cal S}$; (c) $\o(x)=a^{-1}xa$ provided both $x,\, x^{-1}
\in {\Cal S}$; (d) $\o(x)=x$ provided both $x,\, x^{-1} \notin
{\Cal S}$.
\endroster

If $\o$ is of type (W2), then we write $\o=\o({\Cal S}, a)$. It is
obvious that $\a_1$ in the statement of the lemma above is of type (W2).

Combining the conclusions of Steps A--C enables us to immediately
obtain the conclusion of Step ~D.

In Steps E--F, we use, as our main tool, the Whitehead graphs
$\varPhi (Y)$ and $\varPhi_{x_i} (Y)$ of a reduced non-cyclic
word $Y \in F_n$. The standard Whitehead graph
$\varPhi (Y)$ is constructed as follows: Take the vertex set as
${\frak Y}^{\pm 1}=\{\text {letters occurring in} \ Y\}^{\pm 1}$, and connect two vertices $a, \, b \in {\frak Y}^{\pm 1}$ by a non-oriented edge if there is a subword
$ab^{-1}$ or $ba^{-1}$ of (non-cyclic) $Y$. We let $C(a, Y)$
denote the connected component of $\varPhi (Y)$ containing a
vertex $a \in \varPhi(Y)$.

In a similar way, we construct the graph $\varPhi_{x_i} (Y)$,
introduced by Ivanov [2], paying special attention to powers of
$x_i$ in $Y$ as follows: Take the vertex set as ${\frak Y}^{\pm
1} \setminus \{x_i^{\pm 1}\}$, and connect two distinct vertices
$a,\, b \in {\frak Y}^{\pm 1} \setminus \{x_i^{\pm 1}\}$ by a
non-oriented edge if (non-cyclic) $Y$ contains a subword of the
form $ax_i^{\ell}b^{-1}$ or $bx_i^{\ell}a^{-1}$ for some $\ell$.
We will refer to this graph as the generalized Whitehead graph of $(Y, x_i)$. By $C_{x_i}(a, Y)$ we denote the connected component of
$\varPhi_{x_i} (Y)$ containing a vertex $a \in \varPhi_{x_i} (Y)$.

The idea and the techniques used in [2] are developed further in
the present paper.

\heading 2. Proof of the Theorem
\endheading

Thanks to Ivanov's result, it suffices to show that $\<\phi(x_1),
\, \phi(x_2)\>$ is a free factor of rank 2 of $F_n$. Since
$\phi(x_1)$ is primitive in $F_n$, we can assume $\phi(x_1) =
x_1$. Consider the elements $\phi(x_2),\, \phi(x_3)$. Let us
define $W$ as a word in $F_n$ that has the minimum length over
all words of the form $P_1 \, \a \bigl(\phi(x_2)\bigl) \,P_2$, where
$P_1, P_2 \in \<x_1\>$, $\a \in \text {Aut}\, F_n$ and $\a(x_1) =
x_1$. We may assume that $\phi(x_3)$ has the minimum length over
all words of the form
$$Q_1\,\b\bigl(\phi(x_3)\bigl)\,Q_2, 
\tag 1
$$
where $Q_1, \, Q_2 \in \<x_1\>$, $\b
\in \text {Aut} \, F_n$, $\b(x_1) = x_1$, and, for some $S_1,\, S_2
\in \<x_1\>$, $\|S_1\b(W)S_2\| = |W|$.

Let $\g \in \text {Aut} \, F_n$ be used in assuming $\phi(x_3)$
in $(1)$. For such $\g$, there exist $T_1, \, T_2 \in \<x_1\>$
such that $\|T_1\g(W)T_2\| = |W|$. By writing $X \equiv Y$ we mean
the graphical equality (letter-by-letter) of words $X, \, Y \in
F_n$ (whereas $X=Y$ means the equality of $X, \, Y$ in $F_n$). We
then may assume that $\phi(x_2) \equiv \overline{T_1\g(W)T_2}$.
Now put
$$U \equiv \phi(x_2) \quad \text{and} \quad V \equiv \phi(x_3).$$

If $|U| = 1$, then the assertion is obvious. Supposing $|U|>1$, we
shall derive a contradiction. If $|V|=1$, then a contradiction
follows immediately; for if $|V|=1$ then $\<V=\phi(x_3),\,
x_1=\phi(x_1)\>$ would be a free factor of rank 2 of $F_n$, hence
by Ivanov's result, $\phi$ would be an automorphism, but then
$|U|$ would have to be 1. So let $|V|>1$.

Similar arguments to those in [2, Lemmas 1--3] show the existence
of the Whitehead automorphisms $\g_1 = \g_1({\Cal
R}_1,x_1^{-1}),\ \g_2 = \g_2({\Cal R}_2,x_1^{-1})$ such that
$$
\g_1(U)=Ux_1^{-1} \quad \text {and} \quad \g_2(U)= x_1U.
\tag 2
$$
This, in particular, implies that $U$ is a cyclically reduced
word; for otherwise the last letter of $\g_1(U)$ being $x_1^{-1}$
would force the first letter of $\g_1(U)$ to be $x_1$, contrary to
$\g_1(U)=Ux_1^{-1}$. So the word
$$W_1 \equiv (Ux_1^r)^sU^3x_1^rVx_1^rU^{-3}(x_1^{-r}U)^s$$
which we consider in Step A is reduced.

For $y \in {\frak X}^{\pm 1}$, $\t_y$ denotes the inner automorphism of $F_n$ induced by $y$, i.e., $\t_y(x)= yxy^{-1}$ for all $x \in {\frak X}$, while $\pi_y$ denotes an endomorphism of $F_n$ such that $\pi_y(y)=1$ and $\pi_y(x)=x$ for all $x \in {\frak X}^{\pm 1} \setminus \{y^{\pm 1}\}$. 

\medskip
\subheading {Step A} Consider the word
$$W_1 \equiv (Ux_1^r)^sU^3x_1^rVx_1^rU^{-3}(x_1^{-r}U)^s,$$
where $r,\,s$ are integers with $s>|V|+6|U|+10, \
r>|V|+(2s+6)|U|+4s+8$. Since $\phi$ is primitivity preserving, the
word $W_1$ is primitive in $F_n$. Then by the Lemma stated in the
Introduction, there is a Whitehead automorphism $\o_1$ of $F_n$
such that
$$\|\o_1(W_1)\|<|W_1|.
\tag 3
$$
Clearly $\o_1$ has form $(W2)$ , so let $\o_1 = \o_1({\Cal
S}_1,a)$.

\proclaim {Claim 1} $\o_1(x_1)=x_1$.
\endproclaim

\demo {Proof of Claim 1} If $\o_1(x_1) = x_1a$ or $a^{-1}x_1$,
then there are $r$ new occurrences of $a$ in $x_1^r$, hence we
have (for details of the following computation, see [2, p. 95]):
$$
\split \|\o_1(W_1)\| &\ge
s|U|+s(|x_1^r|+r)+3|U|+(|x_1^r|+r)+|V|+(|x_1^r|+r)+3|U|
\\
&\qquad +s(|x_1^{-r}|+r)+s|U|-\{2|V|+2(2s+6)|U|+2(4s+8)\} \\
&> |W_1|,
\endsplit
\tag 4
$$
by the choice of $r$, contradicting $(3)$.

Now let $\o_1(x_1) = a^{-1}x_1a$. Then
$$\|\o_1(W_1)\|=
\|\bigl(\o_1(U)a^{-1}x_1^ra\bigl)^s\o_1(U)^3a^{-1}x_1^ra\,\o_1(V)\,a^{-1}x_1^ra\o_1(U)^{-3}\bigl(a^{-1}x_1^{-r}a\o_1(U)\bigl)^s\|.$$

{\it Case} (i): $\overline{\o_1(U)}$ neither begins with $a^{-1}$
nor ends with $a$. By the
choice of $U$, we have $\|\o_1(U)\| \ge |U|-1$; for if
$\|\o_1(U)\| < |U|-2$, then the automorphism $\t_a\o_1$ of $F_n$
fixes $x_1$ and $\|\t_a\o_1(U)\| < |U|$, contradicting the choice
of $U$; besides, if $\|\o_1(U)\| = |U|-2$, then since
$|a\overline {\o_1(U)}a^{-1}|=|U|$, the word
$a\overline{\o_1(U)}a^{-1}$ also has Whitehead automorphisms of
$F_n$ with property $(2)$, but this contradicts the fact that the
word $a\overline{\o_1(U)}a^{-1}$ is not cyclically reduced. Hence
we have
$$\|\bigl(\o_1(U)a^{-1}x_1^ra\bigl)^s\| \ge |(Ux_1^r)^s|+s \quad \text{and}
\quad \|\bigl(a^{-1}x_1^{-r}a\o_1(U)\bigl)^s\| \ge
|(x_1^{-r}U)^s|+s;$$ thus
$$\split
\|\o_1(W_1)\|
&\ge \{|(Ux_1^r)^s|+s\}+3|U|+|x_1^r|+|V|+|x_1^r|+3|U|+\{|(x_1^{-r}U)^s|+s\}
\\
&\qquad -\{2|V|+12|U|+20\} \\
&> |W_1|,
\endsplit
\tag 5
$$
by the choice of $s$, contradicting $(3)$.

{\it Case} (ii): $\overline{\o_1(U)}$ begins with $a^{-1}$ but does not end
with $a$, say, $\overline{\o_1(U)} \equiv a^{-1}U'$ (note that $U'$
cannot begin with $x_1^{\pm 1}$). By the choice of $U$, $|U'| \ge
|U|-1$; for if $|U'| < |U|-1$ then the automorphism $\t_a\o_1$
fixes $x_1$ and $\|\t_a\o_1(U)\|=|U'a^{-1}| < |U|$, contrary to
the choice of $U$. Moreover, if $|U'| > |U|-1$, then
$\|\bigl(\o_1(U)a^{-1}x_1^ra\bigl)^s\| \ge |(Ux_1^r)^s|+s$ and
$\|\bigl(a^{-1}x_1^{-r}a\o_1(U)\bigl)^s\| \ge |(x_1^{-r}U)^s|+s$;
hence as in $(5)$, we will get a contradiction to $(3)$. So let
$|U'|= |U|-1$. Then the automorphism $\t_a\o_1$ of $F_n$ fixes
$x_1$ and $\|\t_a\o_1(U)\| = |U|$. Hence by the choice of $V$,
$\|\t_a\o_1(V)\| \ge |V|$, i.e., $\|a\,\o_1(V)\,a^{-1}\| \ge |V|$.
But then, from the following observation
$$\split
\|\o_1(W_1)\| &= \|\o_1\bigl((Ux_1^r)^sU^3x_1^r\bigl)\|+
\|a\,\o_1(V)\,a^{-1}\|+ \|\o_1\bigl(x_1^rU^{-3}(x_1^{-r}U)^s\bigl)\|
-2 \\
& \ge \|\o_1\bigl((Ux_1^r)^sU^3x_1^r\bigl)\|+ |V|+
\|\o_1\bigl(x_1^rU^{-3}(x_1^{-r}U)^s\bigl)\|-2,
\endsplit
\tag 6
$$
we have $\|\o_1(W_1)\| \ge |W_1|$, contrary to $(3)$, because
$$\split
\|\o_1\bigl((Ux_1^r)^sU^3x_1^r\bigl)\|
&=\|(a^{-1}U'a^{-1}x_1^ra)^s(a^{-1}U')^3a^{-1}x_1^ra\|
\\
&= |(Ux_1^r)^sU^3x_1^r|+(2s+2)-\{2(s-1)+2\} \\
&= |(Ux_1^r)^sU^3x_1^r|+2,
\endsplit
$$
$$\split
\|\o_1\bigl(x_1^rU^{-3}(x_1^{-r}U)^s\bigl)\|
&=\|a^{-1}x_1^ra({U'}^{-1}a)^3(a^{-1}x_1^{-r}aa^{-1}U')^s\|
\\
&=|x_1^rU^{-3}(x_1^{-r}U)^s|+(2+2s)-(2+2s) \\
&=|x_1^rU^{-3}(x_1^{-r}U)^s|.
\endsplit
$$

{\it Case} (iii): $\overline{\o_1(U)}$ ends with $a$ but does not begin with
$a^{-1}$, say, $\overline{\o_1(U)} \equiv U'a$ (note that $U'$
cannot end with $x_1^{\pm 1}$). The same reason as above enables
us to assume $|U'|=|U|-1$. We then have:
$$\split
\|\o_1\bigl((Ux_1^r)^sU^3x_1^r\bigl)\|
&=\|(U'aa^{-1}x_1^ra)^s(U'a)^3a^{-1}x_1^ra\|
\\
&= |(Ux_1^r)^sU^3x_1^r|+(2s+2)-(2s+2) \\
&= |(Ux_1^r)^sU^3x_1^r|;
\endsplit
$$
$$\split
\|\o_1\bigl(x_1^rU^{-3}(x_1^{-r}U)^s\bigl)\|
&=\|a^{-1}x_1^ra(a^{-1}{U'}^{-1})^3(a^{-1}x_1^{-r}aU'a)^s\|
\\
&= |x_1^rU^{-3}(x_1^{-r}U)^s|+(2+2s)-\{2+2(s-1)\} \\
&= |x_1^rU^{-3}(x_1^{-r}U)^s|+2.
\endsplit
$$
Then reasoning as above, we reach a contradiction to $(3)$.

{\it Case} (iv): $\overline{\o_1(U)}$ both begins with $a^{-1}$ and ends
with $a$, say, $\overline{\o_1(U)} \equiv a^{-1}U'a$ (note that
$U'$ is cyclically reduced and can neither begin nor end with
$x_1^{\pm 1}$). By the choice of $U$, $|U'| \ge |U|$. However, if
$|U'|
> |U|$ then $\|\bigl(\o_1(U)a^{-1}x_1^ra\bigl)^s\| \ge
|(Ux_1^r)^s|+s$ and $\|\bigl(a^{-1}x_1^{-r}a\o_1(U)\bigl)^s\| \ge
 |(x_1^{-r}U)^s|+s$, so as in $(5)$, we will arrive at a contradiction to $(3)$. So let
 $|U'|=|U|$. Then we have:
$$\split
\|\o_1\bigl((Ux_1^r)^sU^3x_1^r\bigl)\|
&=|a^{-1}(U'x_1^r)^s{U'}^3x_1^ra| \\
&=|(Ux_1^r)^sU^3x_1^r|+2; \\
\endsplit
$$
$$\split
\|\o_1\bigl(x_1^rU^{-3}(x_1^{-r}U)^s\bigl)\|
&=|a^{-1}x_1^r{U'}^{-3}(x_1^{-r}U')^sa| \\
&=|x_1^rU^{-3}(x_1^{-r}U)^s|+2.
\endsplit
$$
Moreover, by the choice of $V$, $\|a\o_1(V)a^{-1}\| \ge |V|$. It
then follows from $(6)$ that $\|\o_1(W_1)\| \ge |W_1|$,
contradicting $(3)$. This completes the proof of Claim 1. \qed
\enddemo

By Claim 1 together with the choice of $U$, we have $\|\o_1(U)\|
\ge |U|$.

\proclaim {Claim 2} $\|\o_1(U)\|=|U|$.
\endproclaim

\demo {Proof of Claim 2} By way of contradiction, suppose that
$\|\o_1(U)\|>|U|$. If $\|\bigl(\o_1(U)x_1^r\bigl)^s\|+
\|\bigl(x_1^{-r}\o_1(U)\bigl)^s\|>
|(Ux_1^r)^s|+|(x_1^{-r}U)^s|+2s$, then as in $(5)$, we will have a
contradiction to $(3)$. In order to avoid this contradiction,
$\overline{\o_1(U)}$ must have one of the forms $U'x_1^{\pm 1}$,
$x_1^{\pm 1}U'$ or $x_1^{\pm 1}U'x_1^{\mp 1}$, where $|U'|=|U|$
and $U'$ neither begins nor ends with $x_1^{\pm 1}$; in any case,
$\|x_1\o_1(V)x_1\| \ge |V|+1$ by the choice of $V$ together with
the fact that there can be at most one cancellation in
$x_1\overline{\o_1(V)}x_1$. Now observe:
$$\split
\|\o_1(W_1)\| &= \|\o_1\bigl((Ux_1^r)^sU^3x_1^r\bigl)\|+
\|x_1\o_1(V)x_1\|+ \|\o_1\bigl(x_1^rU^{-3}(x_1^{-r}U)^s\bigl)\|
-2 \\
& \ge \|\o_1\bigl((Ux_1^r)^sU^3x_1^r\bigl)\|+ |V|+
\|\o_1\bigl(x_1^rU^{-3}(x_1^{-r}U)^s\bigl)\|-1.
\endsplit
\tag 7
$$

{\it Case} (i): $\overline{\o_1(U)} \equiv U'x_1$ with $|U'|=|U|$; then
$$\split
\|\o_1\bigl((Ux_1^r)^sU^3x_1^r\bigl)\|
&=\|(U'x_1x_1^r)^s(U'x_1)^3x_1^r\|
\\
&= |(Ux_1^r)^sU^3x_1^r|+(s+3),
\endsplit
$$
$$\split
\|\o_1\bigl(x_1^rU^{-3}(x_1^{-r}U)^s\bigl)\|
&=\|x_1^r(x_1^{-1}{U'}^{-1})^3(x_1^{-r}U'x_1)^s\|
\\
&= |x_1^rU^{-3}(x_1^{-r}U)^s|+(3+s)-\{2+2(s-1)\} \\
&= |x_1^rU^{-3}(x_1^{-r}U)^s|-(s-3),
\endsplit
$$
which yields by $(7)$ that $\|\o_1(W_1)\| \ge |W_1|$,
contradicting $(3)$.

{\it Case} (ii): $\overline{\o_1(U)} \equiv U'x_1^{-1}$ with
$|U'|=|U|$; then
$$\split
\|\o_1\bigl((Ux_1^r)^sU^3x_1^r\bigl)\|
&=\|(U'x_1^{-1}x_1^r)^s(U'x_1^{-1})^3x_1^r\|
\\
&= |(Ux_1^r)^sU^3x_1^r|+(s+3)-(2s+2) \\
&= |(Ux_1^r)^sU^3x_1^r|-(s-1),
\endsplit
$$
$$\split
\|\o_1\bigl(x_1^rU^{-3}(x_1^{-r}U)^s\bigl)\|
&=\|x_1^r(x_1{U'}^{-1})^3(x_1^{-r}U'x_1^{-1})^s\|
\\
&= |x_1^rU^{-3}(x_1^{-r}U)^s|+(3+s);
\endsplit
$$
thus $\|\o_1(W_1)\| \ge |W_1|$ by $(7)$, contradicting $(3)$ as
well.

{\it Case} (iii): $\overline{\o_1(U)} \equiv x_1U'$ with $|U'|=|U|$; then
$$\split
\|\o_1\bigl((Ux_1^r)^sU^3x_1^r\bigl)\|
&=\|(x_1U'x_1^r)^s(x_1U')^3x_1^r\|
\\
&= |(Ux_1^r)^sU^3x_1^r|+(s+3),
\endsplit
$$
$$\split
\|\o_1\bigl(x_1^rU^{-3}(x_1^{-r}U)^s\bigl)\|
&=\|x_1^r({U'}^{-1}x_1^{-1})^3(x_1^{-r}x_1U')^s\|
\\
&= |x_1^rU^{-3}(x_1^{-r}U)^s|+(3+s)-2s \\
&= |x_1^rU^{-3}(x_1^{-r}U)^s|-(s-3);
\endsplit
$$
hence $\|\o_1(W_1)\| \ge |W_1|$ again by $(7)$, contradicting
$(3)$.

{\it Case} (iv): $\overline{\o_1(U)} \equiv x_1^{-1}U'$ with
$|U'|=|U|$; then
$$\split
\|\o_1\bigl((Ux_1^r)^sU^3x_1^r\bigl)\|
&=\|(x_1^{-1}U'x_1^r)^s(x_1^{-1}U')^3x_1^r\|
\\
&= |(Ux_1^r)^sU^3x_1^r|+(s+3)-\{2(s-1)+2\} \\
&= |(Ux_1^r)^sU^3x_1^r|-(s-3),
\endsplit
$$
$$\split
\|\o_1\bigl(x_1^rU^{-3}(x_1^{-r}U)^s\bigl)\|
&=\|x_1^r({U'}^{-1}x_1)^3(x_1^{-r}x_1^{-1}U')^s\|
\\
&= |x_1^rU^{-3}(x_1^{-r}U)^s|+(3+s)-2 \\
&= |x_1^rU^{-3}(x_1^{-r}U)^s|+(s+1);
\endsplit
$$
thus $\|\o_1(W_1)\| \ge |W_1|$ by $(7)$, contrary to $(3)$.

{\it Case} (v): $\overline{\o_1(U)} \equiv x_1U'x_1^{-1}$ with
$|U'|=|U|$; then
$$\split
\|\o_1\bigl((Ux_1^r)^sU^3x_1^r\bigl)\|
&=\|(x_1U'x_1^{-1}x_1^r)^s(x_1U'x_1^{-1})^3x_1^r\|
\\
&= |(Ux_1^r)^sU^3x_1^r|+(2s+6)-(2s+4+2) \\
&= |(Ux_1^r)^sU^3x_1^r|,
\endsplit
$$
$$\split
\|\o_1\bigl(x_1^rU^{-3}(x_1^{-r}U)^s\bigl)\|
&=\|x_1^r(x_1{U'}^{-1}x_1^{-1})^3(x_1^{-r}x_1U'x_1^{-1})^s\|
\\
&= |x_1^rU^{-3}(x_1^{-r}U)^s|+(6+2s)-(4+2s) \\
&= |x_1^rU^{-3}(x_1^{-r}U)^s|+2;
\endsplit
$$
hence $\|\o_1(W_1)\| \ge |W_1|$ by $(7)$, contradicting $(3)$.

{\it Case} (vi): $\overline{\o_1(U)} \equiv x_1^{-1}U'x_1$ with
$|U'|=|U|$; then
$$\split
\|\o_1\bigl((Ux_1^r)^sU^3x_1^r\bigl)\|
&=\|(x_1^{-1}U'x_1x_1^r)^s(x_1^{-1}U'x_1)^3x_1^r\|
\\
&= |(Ux_1^r)^sU^3x_1^r|+(2s+6)-\{2(s-1)+2+4\} \\
&= |(Ux_1^r)^sU^3x_1^r|+2,
\endsplit
$$
$$\split
\|\o_1\bigl(x_1^rU^{-3}(x_1^{-r}U)^s\bigl)\|
&=\|x_1^r(x_1^{-1}{U'}^{-1}x_1)^3(x_1^{-r}x_1^{-1}U'x_1)^s\|
\\
&= |x_1^rU^{-3}(x_1^{-r}U)^s|+(6+2s)-\{2+4+2+2(s-1)\} \\
&= |x_1^rU^{-3}(x_1^{-r}U)^s|;
\endsplit
$$
thus $(7)$ implies that $\|\o_1(W_1)\| \ge |W_1|$, contrary to
$(3)$.

The proof of Claim 2 is complete. \qed
\enddemo

Then by the choice of $U$, $\overline{\o_1(U)}$ neither begins
nor ends with $x_1^{\pm 1}$; besides, $\overline{\o_1(U)}$ has
Whitehead automorphisms with property $(2)$, implying that
$\overline{\o_1(U)}$ is cyclically reduced. Also, by the choice
of $V$, $\|\o_1(V)\| \ge |V|$. In particular, if
$\overline{\o_1(V)}$ either begins or ends with $x_1^{\pm 1}$,
then $\|\o_1(V)\| \ge |V|+1$; if $\overline{\o_1(V)}$ both begins
and ends with $x_1^{\pm 1}$, then $\|\o_1(V)\| \ge |V|+2$.

Now, since $\|\o_1(W_1)\|<|W_1|$, there must be cancellations
in the product $x_1\overline{\o_1(V)}x_1$ in order that
$\|x_1\overline{\o_1(V)}x_1\| < |x_1Vx_1|$. This is possible only
when either $\overline{\o_1(V)} \equiv V_1x_1^{-1}$ or
$\overline{\o_1(V)} \equiv x_1^{-1}V_1$, where $|V_1|=|V|$.
Hence, as our conclusion of Step A, we have the following two
possibilities:
$$
\alignat 3 {\text {\bf (A1)}}\ \o_1(x_1)&=x_1, & \quad
\overline{\o_1(U)} &\equiv
U_1, & \quad \overline{\o_1(V)} &\equiv V_1x_1^{-1}; \\
{\text {\bf (A2)}}\ \o_1(x_1)&=x_1, & \quad \overline{\o_1(U)}
&\equiv U_1, & \quad \overline{\o_1(V)} &\equiv x_1^{-1}V_1,
\endalignat
$$
where $|U_1|=|U|$ and $|V_1|=|V|$.

If (A1) occurs, then the word $VU$ is reduced; for otherwise the
last letter of $\overline{\o_1(V)}$ being $x_1^{-1}$ would force
the first letter of $\overline{\o_1(U)}$ to be $x_1$. Then in the
next step, consider the reduced word
$$W_2 \equiv (U^{-1}x_1^{-r})^sU^{-3}x_1^rVU^3x_1^{-r}(Ux_1^{-r})^s,$$
where $r,\,s$ are integers with $s>|V|+6|U|+10, \
r>|V|+(2s+6)|U|+4s+8$.

On the other hand, if (A2) occurs, then the word $UV$ is reduced;
then in the next step, consider the reduced word
$$W_2' \equiv (x_1^{-r}U)^sx_1^{-r}U^3Vx_1^rU^{-3}(x_1^{-r}U^{-1})^s,$$
where $r,\,s$ are integers with $s>|V|+6|U|+10, \
r>|V|+(2s+6)|U|+4s+8$.

\medskip
\subheading {Step B} Suppose that (A1) occurs. We consider the
word
$$W_2 \equiv (U^{-1}x_1^{-r})^sU^{-3}x_1^rVU^3x_1^{-r}(Ux_1^{-r})^s,$$
where $r,\,s$ are integers with $s>|V|+6|U|+10, \
r>|V|+(2s+6)|U|+4s+8$. Since $\phi$ is primitivity preserving, the
word $W_2$ is primitive in $F_n$; so by the Lemma, there is a
Whitehead automorphism $\o_2 = \o_2({\Cal S}_2,b) $ of $F_n$ such
that
$$\|\o_2(W_2)\|<|W_2|.
\tag 7
$$

\proclaim {Claim 3} $\o_2(x_1)=x_1$.
\endproclaim

\demo {Proof of Claim 3} For the same reason as in $(4)$, $\o_2(x_1)$
cannot be of the form $x_1b$ or $b^{-1}x_1$. So let $\o_2(x_1) =
b^{-1}x_1b$. We have
$$\|\o_2(W_2)\|=
\|\bigl(\o_2(U)^{-1}b^{-1}x_1^{-r}b\bigl)^s\o_2(U)^{-3}b^{-1}x_1^rb\,\o_2(V)\,\o_2(U)^3b^{-1}x_1^{-r}b
\bigl(\o_2(U)b^{-1}x_1^{-r}b\bigl)^s\|.$$ 

{\it Case} (i): $\overline{\o_2(U)}$ neither begins with $b^{-1}$ nor ends with
$b$. Then by the same argument as with $W_1$ of this case, we
arrive at a contradiction to $(7)$.

{\it Case} (ii): $\overline{\o_2(U)}$ begins with $b^{-1}$ but does not end
with $b$, say, $\overline{\o_2(U)} \equiv b^{-1}U'$ ($U'$ cannot
begin with $x_1^{\pm 1}$). Reasoning in the same way as with $W_1$
of this case enables us to assume $|U'|=|U|-1$. Then
$\|b\,\o_2(V)\,b^{-1}\| = \|\t_b\o_2(V)\| \ge |V|$ by the choice of
$V$; moreover
$$\split
\|\o_2\bigl((U^{-1}x_1^{-r})^sU^{-3}x_1^r\bigl)\|
&=\|({U'}^{-1}bb^{-1}x_1^{-r}b)^s({U'}^{-1}b)^3b^{-1}x_1^rb\| \\
&=|(U^{-1}x_1^{-r})^sU^{-3}x_1^r|+(2s+2)-(2s+2) \\
&=|(U^{-1}x_1^{-r})^sU^{-3}x_1^r|,
\endsplit
$$
$$\split
\|\o_2\bigl(U^3x_1^{-r}(Ux_1^{-r})^s\bigl)\|
&=\|(b^{-1}U')^3b^{-1}x_1^{-r}b(b^{-1}U'b^{-1}x_1^{-r}b)^s\| \\
&=|U^3x_1^{-r}(Ux_1^{-r})^s|+(2+2s)-\{2+2(s-1)\} \\
&=|U^3x_1^{-r}(Ux_1^{-r})^s|+2.
\endsplit
$$
It then follows from an observation similar to $(6)$ that
$\|\o_2(W_2)\| \ge |W_2|$, contradicting $(7)$.

{\it Case} (iii): $\overline{\o_2(U)}$ ends with $b$ but does not begin with
$b^{-1}$, say, $\overline{\o_2(U)} \equiv U'b$ ($U'$ cannot end
with $x_1^{\pm 1}$). In this case as well, we may assume
$|U'|=|U|-1$. Then by the choice of $V$, $\|\o_2(V)\| \ge |V|-2$.
In particular, if $\overline{\o_2(V)}$ begins with $b^{-1}$, then
$\|\o_2(V)\| \ge \|\t_b\o_2(V)\| \ge |V|$. So $\|b\, \o_2(V)\| \ge
|V|-1$. Since all cancellations in $\o_2(W_2)$ have the form
$b^{\pm 1}b^{\mp 1} \rightarrow 1$, there can be no cancellation
in $\overline{\o_2(V)}U'b$. In addition, we have
$$\split
\|\o_2\bigl((U^{-1}x_1^{-r})^sU^{-3}x_1^r\bigl)\|
&=\|(b^{-1}{U'}^{-1}b^{-1}x_1^{-r}b)^s(b^{-1}{U'}^{-1})^3b^{-1}x_1^rb\| \\
&=|(U^{-1}x_1^{-r})^sU^{-3}x_1^r|+(2s+2)-\{2(s-1)+2\} \\
&=|(U^{-1}x_1^{-r})^sU^{-3}x_1^r|+2,
\endsplit
$$
$$\split
\|\o_2\bigl(U^3x_1^{-r}(Ux_1^{-r})^s\bigl)\|
&=\|(U'b)^3b^{-1}x_1^{-r}b(U'bb^{-1}x_1^{-r}b)^s\| \\
&=|U^3x_1^{-r}(Ux_1^{-r})^s|+(2+2s)-(2+2s) \\
&=|U^3x_1^{-r}(Ux_1^{-r})^s|.
\endsplit
$$
Therefore,
$$\split
\|\o_2(W_2)\| &=
\|\o_2\bigl((U^{-1}x_1^{-r})^sU^{-3}x_1^r\bigl)\|+ \|b\o_2(V)\|+
\|\o_2\bigl(U^3x_1^{-r}(Ux_1^{-r})^s\bigl)\|
- 1 \\
& \ge \|\o_2\bigl((U^{-1}x_1^{-r})^sU^{-3}x_1^r\bigl)\|+ |V|+
\|\o_2\bigl(U^3x_1^{-r}(Ux_1^{-r})^s\bigl)\|-2 \\
&  \ge |W_2|,
\endsplit
$$
contradicting $(7)$.

{\it Case} (iv): $\overline{\o_2(U)}$ both begins with $b^{-1}$ and ends
with $b$, say, $\overline{\o_2(U)} \equiv b^{-1}U'b$ ($U'$ is
cyclically reduced and can neither begin nor end with $x_1^{\pm
1}$). As with $W_1$ of this case, we arrive at a contradiction to
$(7)$. This completes the proof of Claim 3. \qed
\enddemo

Now, Claim 3 together with the choice of $U$ implies $\|\o_2(U)\|
\ge |U|$.

\proclaim {Claim 4} $\|\o_2(U)\|=|U|$.
\endproclaim

\demo {Proof of Claim 4} By way of contradiction, suppose that
$\|\o_2(U)\|>|U|$. Reasoning in the same way as with $W_1$, $\overline{\o_2(U)}$
must have one of the forms $U'x_1^{\pm 1}$, $x_1^{\pm 1}U'$ or
$x_1^{\pm 1}U'x_1^{\mp 1}$, where $|U'|=|U|$ and $U'$ neither
begins nor ends with $x_1^{\pm 1}$; in any case, we see from the
choice of $V$ that
$$\|\o_2(W_2)\| \ge
\|\o_2\bigl((U^{-1}x_1^{-r})^sU^{-3}x_1^r\bigl)\|+ |V|+
\|\o_2\bigl(U^3x_1^{-r}(Ux_1^{-r})^s\bigl)\|-2. \tag 8
$$

{\it Case} (i): $\overline{\o_2(U)} \equiv U'x_1$ with $|U'|=|U|$;
then
$$\split
\|\o_2\bigl((U^{-1}x_1^{-r})^sU^{-3}x_1^r\bigl)\|
&=\|(x_1^{-1}{U'}^{-1}x_1^{-r})^s(x_1^{-1}{U'}^{-1})^3x_1^r\| \\
&=|(U^{-1}x_1^{-r})^sU^{-3}x_1^r|+(s+3),
\endsplit
$$
$$\split
\|\o_2\bigl(U^3x_1^{-r}(Ux_1^{-r})^s\bigl)\|
&=\|(U'x_1)^3x_1^{-r}(U'x_1x_1^{-r})^s\| \\
&=|U^3x_1^{-r}(Ux_1^{-r})^s|+(3+s)-(2+2s) \\
&=|U^3x_1^{-r}(Ux_1^{-r})^s|-(s-1).
\endsplit
$$
This yields by $(8)$ that $\|\o_2(W_2)\| \ge |W_2|$, contradicting
$(7)$.

{\it Case} (ii): $\overline{\o_2(U)} \equiv U'x_1^{-1}$ with
$|U'|=|U|$; then
$$\split
\|\o_2\bigl((U^{-1}x_1^{-r})^sU^{-3}x_1^r\bigl)\|
&=\|(x_1{U'}^{-1}x_1^{-r})^s(x_1{U'}^{-1})^3x_1^r\| \\
&=|(U^{-1}x_1^{-r})^sU^{-3}x_1^r|+(s+3)-\{2(s-1)+2\} \\
&=|(U^{-1}x_1^{-r})^sU^{-3}x_1^r|-(s-3),
\endsplit
$$
$$\split
\|\o_2\bigl(U^3x_1^{-r}(Ux_1^{-r})^s\bigl)\|
&=\|(U'x_1^{-1})^3x_1^{-r}(U'x_1^{-1}x_1^{-r})^s\| \\
&=|U^3x_1^{-r}(Ux_1^{-r})^s|+(3+s).
\endsplit
$$
Again by $(8)$ $\|\o_2(W_2)\| \ge |W_2|$, contradicting $(7)$.

{\it Case} (iii): $\overline{\o_2(U)} \equiv x_1U'$ with $|U'|=|U|$; then
$$\split
\|\o_2\bigl((U^{-1}x_1^{-r})^sU^{-3}x_1^r\bigl)\|
&=\|({U'}^{-1}x_1^{-1}x_1^{-r})^s({U'}^{-1}x_1^{-1})^3x_1^r\| \\
&=|(U^{-1}x_1^{-r})^sU^{-3}x_1^r|+(s+3)-2 \\
&=|(U^{-1}x_1^{-r})^sU^{-3}x_1^r|+(s+1),
\endsplit
$$
$$\split
\|\o_2\bigl(U^3x_1^{-r}(Ux_1^{-r})^s\bigl)\|
&=\|(x_1U')^3x_1^{-r}(x_1U'x_1^{-r})^s\| \\
&=|U^3x_1^{-r}(Ux_1^{-r})^s|+(3+s)-\{2+2(s-1)\} \\
&=|U^3x_1^{-r}(Ux_1^{-r})^s|-(s-3),
\endsplit
$$
which yields by $(8)$ that $\|\o_2(W_2)\| \ge |W_2|$, contrary to
$(7)$.

{\it Case} (iv): $\overline{\o_2(U)} \equiv x_1^{-1}U'$ with
$|U'|=|U|$; then
$$\split
\|\o_2\bigl((U^{-1}x_1^{-r})^sU^{-3}x_1^r\bigl)\|
&=\|({U'}^{-1}x_1x_1^{-r})^s({U'}^{-1}x_1)^3x_1^r\| \\
&=|(U^{-1}x_1^{-r})^sU^{-3}x_1^r|+(s+3)-2s, \\
&=|(U^{-1}x_1^{-r})^sU^{-3}x_1^r|-(s-3),
\endsplit
$$
$$\split
\|\o_2\bigl(U^3x_1^{-r}(Ux_1^{-r})^s\bigl)\|
&=\|(x_1^{-1}U')^3x_1^{-r}(x_1^{-1}U'x_1^{-r})^s\| \\
&=|U^3x_1^{-r}(Ux_1^{-r})^s|+(3+s);
\endsplit
$$
thus by $(8)$ $\|\o_2(W_2)\| \ge |W_2|$, contrary to $(7)$.

{\it Case} (v): $\overline{\o_2(U)} \equiv x_1U'x_1^{-1}$ with
$|U'|=|U|$; then
$$\split
\|\o_2\bigl((U^{-1}x_1^{-r})^sU^{-3}x_1^r\bigl)\|
&=\|(x_1{U'}^{-1}x_1^{-1}x_1^{-r})^s(x_1{U'}^{-1}x_1^{-1})^3x_1^r\|
\\
&=|(U^{-1}x_1^{-r})^sU^{-3}x_1^r|+(2s+6)-\{2(s-1)+2+4+2\} \\
&=|(U^{-1}x_1^{-r})^sU^{-3}x_1^r|,
\endsplit
$$
$$\split
\|\o_2\bigl(U^3x_1^{-r}(Ux_1^{-r})^s\bigl)\|
&=\|(x_1U'x_1^{-1})^3x_1^{-r}(x_1U'x_1^{-1}x_1^{-r})^s\| \\
&=|U^3x_1^{-r}(Ux_1^{-r})^s|+(6+2s)-\{4+2+2(s-1)\} \\
&=|U^3x_1^{-r}(Ux_1^{-r})^s|+2;
\endsplit
$$
hence $(8)$ implies that $\|\o_2(W_2)\| \ge |W_2|$, contrary to
$(7)$.

{\it Case} (vi): $\overline{\o_2(U)} \equiv x_1^{-1}U'x_1$ with
$|U'|=|U|$; then
$$\split
\|\o_2\bigl((U^{-1}x_1^{-r})^sU^{-3}x_1^r\bigl)\|
&=\|(x_1^{-1}{U'}^{-1}x_1x_1^{-r})^s(x_1^{-1}{U'}^{-1}x_1)^3x_1^r\|
\\
&=|(U^{-1}x_1^{-r})^sU^{-3}x_1^r|+(2s+6)-(2s+4) \\
&=|(U^{-1}x_1^{-r})^sU^{-3}x_1^r|+2,
\endsplit
$$
$$\split
\|\o_2\bigl(U^3x_1^{-r}(Ux_1^{-r})^s\bigl)\|
&=\|(x_1^{-1}U'x_1)^3x_1^{-r}(x_1^{-1}U'x_1x_1^{-r})^s\| \\
&=|U^3x_1^{-r}(Ux_1^{-r})^s|+(6+2s)-(4+2+2s) \\
&=|U^3x_1^{-r}(Ux_1^{-r})^s|;
\endsplit
$$
thus $\|\o_2(W_2)\| \ge |W_2|$ by $(8)$, contradicting $(7)$.

The proof of Claim 4 is complete. \qed
\enddemo

Then $\overline{\o_2(U)}$ is cyclically reduced and neither begins
nor ends with $x_1^{\pm 1}$. Since $\|\o_2(W_2)\|<|W_2|$, there
must be cancellations in the product
$x_1\overline{\o_2(V)}\,\overline{\o_2(U)}$ in order that
$\|x_1\overline{\o_2(V)}\,\overline{\o_2(U)}\| < |x_1VU|$. Here,
since $\overline{\o_2(U)}$ does not begin with $x_1^{-1}$, there
is precisely one cancellation in the product
$x_1\overline{\o_2(V)}\,\overline{\o_2(U)}$. If the cancellation
occurs in the product $x_1\overline{\o_2(V)}$, then the only
possibility is $\overline{\o_2(V)} \equiv x_1^{-1}V_2$, where
$|V_2|=|V|$; on the other hand, if the cancellation occurs in the
product $\overline{\o_2(V)}\,\overline{\o_2(U)}$, then the only
possibility is that $\overline{\o_2(V)} \equiv V_3a$ and
$\overline{\o_2(U)} \equiv a^{-1}U_3$, where $|V|-1 \le |V_3| \le
|V|\,$, $|U_3|=|U|-1$ and $a \not = x_1^{\pm 1}$.

Similarly repeating Step B for $W_2'$ if (A2) occurs, we conclude
that either $\overline{\o_2'(V)} \equiv V_2'x_1^{-1}$ or
$\overline{\o_2'(V)} \equiv b^{-1}V_3'$, where $\|\o_2'(W_2')\| <
|W_2'|\,$, $|V_2'|=|V|$, $|V|-1 \le |V_3'| \le |V|$ and $b$ is
the last letter of $\overline{\o_2'(U)}$. Consequently Step ~B
gives us the following four possibilities:
$$
\alignat 3 {\text {\bf (B1)}}\ \text {(A1)}, \quad \o_2(x_1)&=x_1,
& \quad
\overline{\o_2(U)} & \equiv U_2, & \quad \overline{\o_2(V)} & \equiv x_1^{-1}V_2; \\
{\text {\bf (B2)}}\ \text {(A1)}, \quad \o_2(x_1)&=x_1, & \quad
\overline{\o_2(U)} & \equiv a^{-1}U_3, & \quad \overline{\o_2(V)} & \equiv V_3a; \\
{\text {\bf (B3)}}\ \text {(A2)}, \quad \o_2'(x_1)&=x_1, & \quad
\overline{\o_2'(U)} & \equiv U_2', & \quad \overline{\o_2'(V)} &
\equiv
V_2'x_1^{-1}; \\
{\text {\bf (B4)}}\ \text {(A2)}, \quad \o_2'(x_1)&=x_1, & \quad
\overline{\o_2'(U)} & \equiv U_3'b, & \quad \overline{\o_2'(V)} &
\equiv
b^{-1}V_3', \\
\endalignat
$$
where $|U_2|=|U_2'|=|U|\,$, $|V_2|=|V_2'|=|V|\,$,
$|U_3|=|U_3'|=|U|-1\,$, $|V|-1 \le |V_3|, |V_3'| \le |V|$ and
$a,\, b \not = x_1^{\pm 1}$.

If (B1) occurs, then the word $UV$ is reduced; so the word $UVU$
is reduced, since the word $VU$ is already reduced by (A1). Also,
if (B3) occurs, then the word $VU$ is reduced; so the word $UVU$
is reduced by (A2). Hence if (B1) or (B3) occurs, then in the
next step, consider the reduced word
$$W_3 \equiv (x_1^rU)^sx_1^rU^3VU^3x_1^{-r}(U^{-1}x_1^r)^s,$$
where $r,\,s$ are integers with $s>|V|+6|U|+10, \
r>|V|+(2s+6)|U|+4s+8$.

\medskip
\subheading {Step C} Suppose that either (B1) or (B3) occurs. We
consider the word
$$W_3 \equiv (x_1^rU)^sx_1^rU^3VU^3x_1^{-r}(U^{-1}x_1^r)^s,$$
where $r,\,s$ are integers with $s>|V|+6|U|+10, \
r>|V|+(2s+6)|U|+4s+8$. The word $W_3$ is primitive in $F_n$; so
by the Lemma, there is a Whitehead automorphism $\o_3= \o_3({\Cal
S}_3,c) $ of $F_n$ such that
$$\|\o_3(W_3)\|<|W_3|.
\tag 9
$$

\proclaim {Claim 5} $\o_3(x_1)=x_1$.
\endproclaim

\demo {Proof of Claim 5} For the same reason as in $(4)$,
$\o_3(x_1)$ cannot be of the form $x_1c$ or $c^{-1}x_1$. Then let
$\o_3(x_1) = c^{-1}x_1c$. We have
$$\|\o_3(W_3)\|=
\|\bigl(c^{-1}x_1^rc\,\o_3(U)\bigl)^sc^{-1}x_1^rc\,\o_3(U)^3\,\o_3(V)\,\o_3(U)^3c^{-1}x_1^{-r}c\,\bigl(\o_3(U)^{-1}c^{-1}x_1^rc\bigl)^s\|.$$

{\it Case} (i): $\overline{\o_3(U)}$ neither
begins with $c^{-1}$ nor ends with $c$; then by the same argument
as with $W_1$ of this case, we get a contradiction to $(9)$.

{\it Case} (ii): $\overline{\o_3(U)}$ begins with $c^{-1}$ but does not end
with $c$, say, $\overline{\o_3(U)} \equiv c^{-1}U'$ ($U'$ cannot
begin with $x_1^{\pm 1}$). Reasoning in the same way as with $W_1$
of this case, we may assume $|U'|=|U|-1$. Then by the choice of
$V$, $\|\o_3(V)\| \ge |V|-2$. In particular, if
$\overline{\o_3(V)}$ ends with $c$, then $\|\o_3(V)\| \ge
\|\t_c\o_3(V)\| \ge |V|$. So $\|\o_3(V)c^{-1}\| \ge |V|-1$. Since
there can be no cancellation in $c^{-1}U'\overline{\o_3(V)}$, we
have:
$$\split
\|\o_3(W_3)\| & = \|\o_3\bigl((x_1^rU)^sx_1^rU^3\bigl)\|+
\|\o_3(V)c^{-1}\|+
\|\o_3\bigl(U^3x_1^{-r}(U^{-1}x_1^r)^s\bigl)\| -1 \\
& \ge \|\o_3\bigl((x_1^rU)^sx_1^rU^3\bigl)\|+ |V|+
\|\o_3\bigl(U^3x_1^{-r}(U^{-1}x_1^r)^s\bigl)\|-2.
\endsplit
\tag 10
$$
Here, we see:
$$\split
\|\o_3\bigl((x_1^rU)^sx_1^rU^3\bigl)\|
&=\|(c^{-1}x_1^rcc^{-1}U')^sc^{-1}x_1^rc(c^{-1}U')^3\| \\
&=|(x_1^rU)^sx_1^rU^3|+(2s+2)-(2s+2) \\
&=|(x_1^rU)^sx_1^rU^3|,
\endsplit
$$
$$\split
\|\o_3\bigl(U^3x_1^{-r}(U^{-1}x_1^r)^s\bigl)\|
&=\|(c^{-1}U')^3c^{-1}x_1^{-r}c({U'}^{-1}cc^{-1}x_1^rc)^s\| \\
&=|U^3x_1^{-r}(U^{-1}x_1^r)^s|+(2+2s)-2s \\
&=|U^3x_1^{-r}(U^{-1}x_1^r)^s|+2,
\endsplit
$$
yielding $\|\o_3(W_3)\| \ge |W_3|$ by $(10)$, contrary to $(9)$.

{\it Case} (iii): $\overline{\o_3(U)}$ ends with $c$ but does not begin with
$c^{-1}$, say, $\overline{\o_3(U)} \equiv U'c$ ($U'$ cannot end
with $x_1^{\pm 1}$). In this case as well, we may assume
$|U'|=|U|-1$. As above, $\|c \, \o_3(V)\| \ge |V|-1$ and there can be
no cancellation in $\o_3(V)U'c$. Moreover, we have:
$$\split
\|\o_3\bigl((x_1^rU)^sx_1^rU^3\bigl)\|
&=\|(c^{-1}x_1^rcU'c)^sc^{-1}x_1^rc(U'c)^3\| \\
&=|(x_1^rU)^sx_1^rU^3|+(2s+2)-\{2(s-1)+2\} \\
&=|(x_1^rU)^sx_1^rU^3|+2,
\endsplit
$$
$$\split
\|\o_3\bigl(U^3x_1^{-r}(U^{-1}x_1^r)^s\bigl)\|
&=\|(U'c)^3c^{-1}x_1^{-r}c(c^{-1}{U'}^{-1}c^{-1}x_1^rc)^s\| \\
&=|U^3x_1^{-r}(U^{-1}x_1^r)^s|+(2+2s)-\{2+2+2(s-1)\} \\
&=|U^3x_1^{-r}(U^{-1}x_1^r)^s|;
\endsplit
$$
so that $\|\o_3(W_3)\| \ge |W_3|$ by an observation similar to
$(10)$, contradicting $(9)$.

{\it Case} (iv): $\overline{\o_3(U)}$ both begins with $c^{-1}$ and ends
with $c$, say, $\overline{\o_3(U)} \equiv c^{-1}U'c$ ($U'$ is
cyclically reduced and can neither begin nor end with $x_1^{\pm
1}$). Then as with $W_1$ of this case, we reach a contradiction to
$(9)$. This completes the proof of Claim ~5. \qed
\enddemo

Claim 5 implies by the choice of $U$ that $\|\o_3(U)\| \ge |U|$.

\proclaim {Claim 6} $\|\o_3(U)\| = |U|$.
\endproclaim

\demo {Proof of Claim 6} Suppose on the contrary that $\|\o_3(U)\|
> |U|$. Reasoning as with $W_1$, $\overline{\o_3(U)}$ must have one
of the forms $U'x_1^{\pm 1}$, $x_1^{\pm 1}U'$ or $x_1^{\pm
1}U'x_1^{\mp 1}$, where $|U'|=|U|$ and $U'$ neither begins nor
ends with $x_1^{\pm 1}$; in any case, it follows from the choice
of $V$ that
$$
\|\o_3(W_3)\| \ge \|\o_3\bigl((x_1^rU)^sx_1^rU^3\bigl)\| + |V| +
\|\o_3\bigl(U^3x_1^{-r}(U^{-1}x_1^r)^s\bigl)\|-2. \tag 11
$$

{\it Case} (i): $\overline{\o_3(U)} \equiv U'x_1$ with $|U'|=|U|$;
then
$$\split
\|\o_3\bigl((x_1^rU)^sx_1^rU^3\bigl)\|
&=\|(x_1^rU'x_1)^sx_1^r(U'x_1)^3\| \\
&=|(x_1^rU)^sx_1^rU^3|+(s+3),
\endsplit
$$
$$\split
\|\o_3\bigl(U^3x_1^{-r}(U^{-1}x_1^r)^s\bigl)\|
&=\|(U'x_1)^3x_1^{-r}(x_1^{-1}{U'}^{-1}x_1^r)^s\| \\
&=|U^3x_1^{-r}(U^{-1}x_1^r)^s|+(3+s)-\{2+2(s-1)\} \\
&=|U^3x_1^{-r}(U^{-1}x_1^r)^s|-(s-3).
\endsplit
$$
This yields by $(11)$ that $\|\o_3(W_3)\| \ge |W_3|$,
contradicting $(9)$.

{\it Case} (ii): $\overline{\o_3(U)} \equiv U'x_1^{-1}$ with
$|U'|=|U|$; then
$$\split
\|\o_3\bigl((x_1^rU)^sx_1^rU^3\bigl)\|
&=\|(x_1^rU'x_1^{-1})^sx_1^r(U'x_1^{-1})^3\| \\
&=|(x_1^rU)^sx_1^rU^3|+(s+3)-\{2(s-1)+2\} \\
&=|(x_1^rU)^sx_1^rU^3|-(s-3),
\endsplit
$$
$$\split
\|\o_3\bigl(U^3x_1^{-r}(U^{-1}x_1^r)^s\bigl)\|
&=\|(U'x_1^{-1})^3x_1^{-r}(x_1{U'}^{-1}x_1^r)^s\| \\
&=|U^3x_1^{-r}(U^{-1}x_1^r)^s|+(3+s)-2 \\
&=|U^3x_1^{-r}(U^{-1}x_1^r)^s|+(s+1);
\endsplit
$$
hence $\|\o_3(W_3)\| \ge |W_3|$ by $(11)$, contradicting $(9)$ as
well.

{\it Case} (iii): $\overline{\o_3(U)} \equiv x_1U'$ with $|U'|=|U|$; then
$$\split
\|\o_3\bigl((x_1^rU)^sx_1^rU^3\bigl)\|
&=\|(x_1^rx_1U')^sx_1^r(x_1U')^3\| \\
&=|(x_1^rU)^sx_1^rU^3|+(s+3),
\endsplit
$$
$$\split
\|\o_3\bigl(U^3x_1^{-r}(U^{-1}x_1^r)^s\bigl)\|
&=\|(x_1U')^3x_1^{-r}({U'}^{-1}x_1^{-1}x_1^r)^s\| \\
&=|U^3x_1^{-r}(U^{-1}x_1^r)^s|+(3+s)-2s \\
&=|U^3x_1^{-r}(U^{-1}x_1^r)^s|-(s-3),
\endsplit
$$
which implies by $(11)$ that $\|\o_3(W_3)\| \ge |W_3|$, contrary
to $(9)$.

{\it Case} (iv): $\overline{\o_3(U)} \equiv x_1^{-1}U'$ with
$|U'|=|U|$; then
$$\split
\|\o_3\bigl((x_1^rU)^sx_1^rU^3\bigl)\|
&=\|(x_1^rx_1^{-1}U')^sx_1^r(x_1^{-1}U')^3\| \\
&=|(x_1^rU)^sx_1^rU^3|+(s+3)-(2s+2) \\
&=|(x_1^rU)^sx_1^rU^3|-(s-1),
\endsplit
$$
$$\split
\|\o_3\bigl(U^3x_1^{-r}(U^{-1}x_1^r)^s\bigl)\|
&=\|(x_1^{-1}U')^3x_1^{-r}({U'}^{-1}x_1x_1^r)^s\| \\
&=|U^3x_1^{-r}(U^{-1}x_1^r)^s|+(3+s);
\endsplit
$$
thus by $(11)$ $\|\o_3(W_3)\| \ge |W_3|$, contrary to $(9)$.

{\it Case} (v): $\overline{\o_3(U)} \equiv x_1U'x_1^{-1}$ with
$|U'|=|U|$; then
$$\split
\|\o_3\bigl((x_1^rU)^sx_1^rU^3\bigl)\|
&=\|(x_1^rx_1U'x_1^{-1})^sx_1^r(x_1U'x_1^{-1})^3\| \\
&=|(x_1^rU)^sx_1^rU^3|+(2s+6)-\{2(s-1)+2+4\} \\
&=|(x_1^rU)^sx_1^rU^3|+2,
\endsplit
$$
$$\split
\|\o_3\bigl(U^3x_1^{-r}(U^{-1}x_1^r)^s\bigl)\|
&=\|(x_1U'x_1^{-1})^3x_1^{-r}(x_1{U'}^{-1}x_1^{-1}x_1^r)^s\| \\
&=|U^3x_1^{-r}(U^{-1}x_1^r)^s|+(6+2s)-(4+2+2s) \\
&=|U^3x_1^{-r}(U^{-1}x_1^r)^s|;
\endsplit
$$
hence by $(11)$ $\|\o_3(W_3)\| \ge |W_3|$, contrary to $(9)$.

{\it Case} (vi): $\overline{\o_3(U)} \equiv x_1^{-1}U'x_1$ with
$|U'|=|U|$; then
$$\split
\|\o_3\bigl((x_1^rU)^sx_1^rU^3\bigl)\|
&=\|(x_1^rx_1^{-1}U'x_1)^sx_1^r(x_1^{-1}U'x_1)^3\| \\
&=|(x_1^rU)^sx_1^rU^3|+(2s+6)-(2s+2+4) \\
&=|(x_1^rU)^sx_1^rU^3|,
\endsplit
$$
$$\split
\|\o_3\bigl(U^3x_1^{-r}(U^{-1}x_1^r)^s\bigl)\|
&=\|(x_1^{-1}U'x_1)^3x_1^{-r}(x_1^{-1}{U'}^{-1}x_1x_1^r)^s\| \\
&=|U^3x_1^{-r}(U^{-1}x_1^r)^s|+(6+2s)-\{4+2+2(s-1)\} \\
&=|U^3x_1^{-r}(U^{-1}x_1^r)^s|+2;
\endsplit
$$
thus $(11)$ implies that $\|\o_3(W_3)\| \ge |W_3|$, contradicting
$(9)$.

The proof of Claim 6 is complete. \qed
\enddemo

Then $\overline{\o_3(U)}$ is cyclically reduced and neither begins
nor ends with $x_1^{\pm 1}$. Since $\|\o_3(W_3)\|<|W_3|$, there
must be cancellations in the product
$\overline{\o_3(U)}\,\overline{\o_3(V)}\,\overline{\o_3(U)}$ so
that
$\|\overline{\o_3(U)}\,\overline{\o_3(V)}\,\overline{\o_3(U)}\| <
|UVU|$. Here, since $\overline{\o_3(U)}$ is cyclically reduced,
there is exactly one cancellation in the product
$\overline{\o_3(U)}\,\overline{\o_3(V)}\,\overline{\o_3(U)}$. If
the cancellation occurs in the product
$\overline{\o_3(U)}\,\overline{\o_3(V)}$, then the only
possibility is that $\overline{\o_3(U)} \equiv U_4c$ and
$\overline{\o_3(V)} \equiv c^{-1}V_4$, where $|U_4|=|U|-1\,$,
$|V|-1 \le |V_4| \le |V|$ and $c \not = x_1^{\pm 1}$; while if the
cancellation occurs in the product
$\overline{\o_3(V)}\,\overline{\o_3(U)}$, then the only
possibility is that $\overline{\o_3(U)} \equiv d^{-1}U_4'$ and
$\overline{\o_3(V)} \equiv V_4'd$, where $|U_4'|=|U|-1\,$, $|V|-1
\le |V_4'| \le |V|$ and $d \not = x_1^{\pm 1}$. Consequently,
Step C provides us with the following four possibilities:
$$
\alignat 3 {\text {\bf (C1)}}\ \text {(B1)},\quad \o_3(x_1)&=x_1,
&\quad
\overline{\o_3(U)} &\equiv U_4c, &\quad \overline{\o_3(V)} &\equiv c^{-1}V_4; \\
{\text {\bf (C2)}}\ \text {(B1)},\quad \o_3(x_1)&=x_1, &\quad
\overline{\o_3(U)} &\equiv d^{-1}U_4', &\quad \overline{\o_3(V)}
&\equiv
V_4'd; \\
{\text {\bf (C3)}}\ \text {(B3)},\quad \o_3(x_1)&=x_1, &\quad
\overline{\o_3(U)} &\equiv U_4c, &\quad \overline{\o_3(V)} &\equiv
c^{-1}V_4; \\
{\text {\bf (C4)}}\ \text {(B3)},\quad \o_3(x_1)&=x_1, &\quad
\overline{\o_3(U)} &\equiv d^{-1}U_4', &\quad \overline{\o_3(V)}
&\equiv V_4'd,
\endalignat
$$
where $|U_4|=|U_4'|=|U|-1\,$, $|V|-1 \le |V_4|, |V_4'| \le |V|$
and $c, \, d \not = x_1^{\pm 1}$.

\medskip
\subheading {Step D} Combining the results of Steps A, B and C
proves the existence of a pair $\bigl(\a=\a({\Cal P},
x_1^{-1}),\, \b=\b({\Cal Q}, e)\bigl)$ or $\bigl(\a'=\a'({\Cal
P}', x_1),\, \b'=\b'({\Cal Q}', e')\bigl)$ of Whitehead
automorphisms of $F_n$ such that
$$
\alignat 3 {\text {\bf (D1)}} &\ \a(x_1)=x_1, &\quad
\overline{\a(U)} &\equiv U_5, &\quad
\overline{\a(V)} &\equiv V_5x_1^{-1},\\
&\ \b(x_1)=x_1, &\quad \overline{\b(U)} &\equiv e^{-1}U_6,
&\quad \overline{\b(V)} &\equiv V_6e;\\
{\text {\bf (D2)}} &\ \a'(x_1)=x_1, &\quad \overline{\a'(U)}
&\equiv
U_5', &\quad \overline{\a'(V)} &\equiv x_1^{-1}V_5',\\
&\ \b'(x_1)=x_1, &\quad \overline{\b'(U)} &\equiv U_6'e', &\quad
\overline{\b'(V)} &\equiv {e'}^{-1}V_6',
\endalignat
$$
where $|U_5|=|U_5'|=|U|\,$, $|V_5|=|V_5'|=|V|\,$,
$|U_6|=|U_6'|=|U|-1\,$, $|V|-1 \le |V_6|, |V_6'| \le |V|$ and $e,
\, e' \not = x_1^{\pm 1}$.

\medskip
\subheading {Step E} Let us write indices as follows: $U_7
\equiv U$, $V_7 \equiv V$, $\a_7=\a$, $\a_7'=\a'$, and
$$
\split
&U_8 \equiv \a_7(U_7) \quad \text{and} \quad V_8 \equiv
\overline{\a_7(V_7)x_1} \quad \text{provided (D1) occurs
for}\ U,\, V; \\
&U_8 \equiv \a_7'(U_7) \quad \text{and} \quad V_8 \equiv
\overline{x_1\a_7'(V_7)} \quad \text{provided (D2) occurs
for}\ U,\, V.
\endsplit
$$
It is easy to see that the words $U_8$ and $V_8$ have the same
properties as $U$ and $V$ do, respectively; hence there exists a pair
$(\a_8, \b_8)$ or $(\a_8', \b_8')$ of Whitehead automorphisms of
$F_n$ that correspond to (D1) or (D2), respectively, relative
to the words $U_8$ and $V_8$. Then let us put
$$
\split
&U_9 \equiv \a_8(U_8) \quad \text{and} \quad V_9 \equiv
\overline{\a_8(V_8)x_1} \quad \text{provided (D1) occurs
for}\ U_8, \, V_8; \\
&U_9 \equiv \a_8'(U_8) \quad \text{and} \quad V_9 \equiv
\overline{x_1\a_8'(V_8)} \quad \text{provided (D2) occurs
for}\ U_8,\, V_8.
\endsplit
$$
Continue this
process until we get the sequences
$$\tilde{\a}_7,\, \tilde{\a}_8, \, \dots, \, \tilde{\a}_l \quad \text
{and} \quad (U_7,V_7), \, (U_8,V_8), \, \dots, \, (U_l,V_l),
$$
where $\tilde{\a}_j = \a_j$ or $\a_j'$ and $l >
6+(2n)^{|U|+|V|}$. We then notice that there is a repetition in
the second sequence, say,
$$(U_{k_1}, V_{k_1})=(U_{k_2}, V_{k_2})$$
with $k_1 < k_2$.

Put $\tilde{\a}(k_2,k_1) =
\tilde{\a}_{k_2-1}\tilde{\a}_{k_2-2}\cdots\tilde{\a}_{k_1}$; then
$$\tilde{\a}(k_2,k_1)(x_1)=x_1, \quad
\tilde{\a}(k_2,k_1)(U_{k_1})=U_{k_1} \quad \text{and} \quad
\tilde{\a}(k_2,k_1)(V_{k_1})=x_1^{-m_1}V_{k_1}x_1^{-m_2}, \tag 12
$$
where $m_1, \, m_2 \ge 0$ and $m_1+m_2=k_2-k_1$. If $m_1=0$, then
$(12)$ implies that there is a Whitehead automorphism
$\d_1=\d_1({\Cal P}_1, x_1^{-1})$ such that
$\d_1(U_{k_1})=U_{k_1}$ and $\d_1(V_{k_1})=V_{k_1}x_1^{-1}$; if
$m_2=0$, then again by $(12)$, there is a Whitehead automorphism
$\d_2=\d_2({\Cal P}_2, x_1)$ such that $\d_2(U_{k_1})=U_{k_1}$ and
$\d_2(V_{k_1})=x_1^{-1}V_{k_1}.$ Now let both $m_1, \, m_2 \not =
0$. In view of $(12)$, there exist two Whitehead automorphisms
$\d_1'=\d_1'({\Cal P}_1', x_1^{-1})$ and $\d_2'=\d_2'({\Cal P}_2',
x_1)$ such that $\d_2'\d_1'(U_{k_1})=U_{k_1}$ and
$\d_2'\d_1'(V_{k_1})=x_1^{-1}V_{k_1}x_1^{-1}$. Here, we may assume
${\Cal P}_1' \cap {\Cal P}_2' = \emptyset$. It then follows that
$\d_1'(U_{k_1})=U_{k_1}$, $\d_2'(U_{k_1})=U_{k_1}$,
$\d_1'(V_{k_1})=V_{k_1}x_1^{-1}$ and
$\d_2'(V_{k_1})=x_1^{-1}V_{k_1}$. Hence, if the words $U_{k_1}$
and $V_{k_1}$ have a pair $(\a_{k_1}, \b_{k_1})$ of Whitehead
automorphisms of type (D1), then we may assume that $\a_{k_1}$
is such that
$$\a_{k_1}(x_1)=x_1, \quad \a_{k_1}(U_{k_1})=U_{k_1} \quad \text{and} \quad
\a_{k_1}(V_{k_1})=V_{k_1}x_1^{-1};$$ if the words $U_{k_1}$ and
$V_{k_1}$ have a pair $(\a_{k_1}', \b_{k_1}')$ of Whitehead
automorphisms of type (D2), then we may assume that $\a_{k_1}'$
is such that
$$\a_{k_1}'(x_1)=x_1, \quad \a_{k_1}'(U_{k_1})=U_{k_1} \quad \text{and} \quad
\a_{k_1}'(V_{k_1})=x_1^{-1}V_{k_1}.$$ Since
$\varPhi_{x_1}(U)=\varPhi_{x_1}(U_{k_1})$ and
$\varPhi_{x_1}(V)=\varPhi_{x_1}(V_{k_1})$, where
$\varPhi_{x_1}(Y)$ is the generalized Whitehead graph of $(Y, x_1)$, we can finally assume that in (D1),
$$
\a(x_1)=x_1,\quad \a(U)=U \quad \text {and} \quad \a(V)=Vx_1^{-1},
$$
and that in (D2),
$$
\a'(x_1)=x_1,\quad \a'(U)=U \quad \text {and} \quad
\a'(V)=x_1^{-1}V.
$$

Consequently, Step E enables us to refine upon the result of Step
D as follows: There exists a pair $\bigl(\a=\a({\Cal P},
x_1^{-1}),\, \b=\b({\Cal Q}, e)\bigl)$ or $\bigl(\a'=\a'({\Cal
P}', x_1),\, \b'=\b'({\Cal Q}', e')\bigl)$ of Whitehead
automorphisms of $F_n$ such that
$$
\alignat 3 {\text {\bf (E1)}} &\ \a(x_1)=x_1, &\quad
\overline{\a(U)} &\equiv U, &\quad
\overline{\a(V)} &\equiv Vx_1^{-1},\\
&\ \b(x_1)=x_1, &\quad \overline{\b(U)} &\equiv e^{-1}U_6,
&\quad \overline{\b(V)} &\equiv V_6e;\\
{\text {\bf (E2)}} &\ \a'(x_1)=x_1, &\quad \overline{\a'(U)}
&\equiv
U, &\quad \overline{\a'(V)} &\equiv x_1^{-1}V,\\
&\ \b'(x_1)=x_1, &\quad \overline{\b'(U)} &\equiv U_6'e', &\quad
\overline{\b'(V)} &\equiv {e'}^{-1}V_6',
\endalignat
$$
where $|U_6|=|U_6'|=|U|-1\,$, $|V|-1 \le |V_6|, |V_6'| \le |V|$
and $e, \, e' \not = x_1^{\pm 1}$.

\medskip
\subheading {Step F} Suppose that the words $U$ and $V$ have a
pair $(\a, \b)$ of Whitehead automorphisms of type (E1) (the
case (E2) is similar). Let $y$ and $z$ be the first and
the last letter of the word $U$, respectively. Recall that
$\b=\b({\Cal Q}, e)$ with $e \not = x_1^{\pm 1}$. Obviously both
$x_1$, $x_1^{-1} \notin {\Cal Q}$, because $\b(x_1)=x_1$. Then
let us define
$$
\split {\Cal T}_1=\{b \in {\Cal Q} \cap
\varPhi(U)\,|\,&\text{there\ exists\ a\ path}\
p_b=b-y^{-1}\ \text{of\ length}\ \ge 0\ \text{in}\ \varPhi(U)\\
&\text{that\ connects}\ b\ \text{with}\ y^{-1}\ \text{and\ does\ not\
pass\ through}\ x_1^{\pm
1}\};\\
{\Cal T}_2=[{\Cal Q} \cap \varPhi(U)] \,\setminus \, {\Cal T}_1,
\endsplit
$$
where $\varPhi(U)$ is the standard Whitehead graph of $U$.

Consider the Whitehead automorphisms
$$\eta_1=\eta_1({\Cal T}_1 \cup \{e\},\,e) \quad \text {and} \quad \eta_2=\eta_2({\Cal
T}_2 \cup \{e\},\,e).$$ If $e \notin {\Cal T}_1$, then it is not
hard to see that $\|\eta_1(U)\| > |U|$. It then follows from the
observation
$$|U|=\|\b(U)\|=\|\eta_1(U)\|+\|\eta_2(U)\|-|U|
\tag 13
$$
that $\|\eta_2(U)\| < |U|$, contrary to the choice of $U$ (note
$\eta_2(x_1)=x_1$ because $x_1^{\pm 1} \notin {\Cal T}_2
\subseteq {\Cal Q}$). So the vertex $e$ has to be inside ${\Cal
T}_1$, i.e.,
$$e \in {\Cal T}_1.$$

\proclaim {Claim 7} ${\Cal T}_2 = \emptyset$.
\endproclaim

\demo {Proof of Claim 7} Suppose on the contrary that ${\Cal T}_2
\not = \emptyset$. Then since $e \in {\Cal T}_1$, i.e., $e \notin
{\Cal T}_2$, $\|\eta_2(U)\| \ge |U|$. Here, if $\|\eta_2(U)\|
> |U|$, then $\|\eta_1(U)\| < |U|$ by $(13)$. This contradicts the choice of $U$. Hence
$\|\eta_2(U)\|=|U|$; this can happen only when
$${\Cal T}_2=\bigcup_{f \in
{\Cal T}_2}C_e(f, U) \quad \text{and} \quad z \notin {\Cal T}_2,$$
where $C_e(f,U)$ is the connected component of $\varPhi_e(U)$
containing $f$.

Here, assume that $C_e(f, U) \not = C(f, U)$ for some $f \in
{\Cal T}_2$, where $C(f, U)$ is the connected component of
$\varPhi(U)$ containing $f$. This is possible only when $e$ or
$e^{-1} \in C(f, U)$. If $e \in C(f, U)$, then there exists a path
$r_f=f-e$ of length $\ge 0$ in $\varPhi(U)$ which connects $f$
with $e$. Since $x_1^{\pm 1} \notin C_e(f, U) \subseteq {\Cal
T}_2 \subseteq {\Cal Q}$, the path $r_f$ does not pass through
$x_1^{\pm 1}$. But then, since $e \in {\Cal T}_1$, such a vertex
$f$ must belong to ${\Cal T}_1$, contrary to $f \in
{\Cal T}_2$. Hence we must have $e \notin C(f, U)$ and
$e^{-1} \in C(f, U)$. This implies that $C_e(f, U)$ contains at
least one of the adjacent vertices of $e$ in $\varPhi(U)$, say, $g$. Of
course $g \not = x_1^{\pm 1}$, because $x_1^{\pm 1} \notin C_e(f,
U) \subseteq {\Cal T}_2 \subseteq {\Cal Q}$. But then, by the
definition of ${\Cal T}_1$, we have $g \in {\Cal T}_1$, which
contradicts $g \in C_e(f, U) \subseteq {\Cal T}_2$.
So it is proved that
$$C_e(f, U)=C(f, U)\quad \text{for\ all}\ f \in {\Cal T}_2.$$

Now choose an arbitrary $C(f, U)= C_e(f, U)$ which is contained in
${\Cal T}_2$. Assume that there is a vertex, say, $b$, in $C(f,U)$
such that $b \in C(f,U)$ but $b^{-1} \notin C(f,U)$. Then for the
Whitehead automorphism $$\e=\e\bigl(C(f, U),\, b),$$ we have
$\e(x_1)=x_1$, because $x_1^{\pm 1} \notin C(f, U) \subseteq {\Cal
T}_2 \subseteq {\Cal Q}$. Furthermore, we see from $y^{-1}, \, z \notin C(f, U) \subseteq {\Cal T}_2$ that
$\|\e(U)\|=\|\pi_b(U)\| < |U|$ (recall that $\pi_b$ removes all $b^{\pm 1}$'s). This contradicts the choice of $U$. Thus if $b \in
C(f,U)$, then $b^{-1}$ also has to occur in $C(f,U)$, which
forces us to have $C(f,U)=\varPhi(U)$, so that $e \in C(f,
U)$, contrary to $e \notin C_e(f, U)=C(f, U)$. This completes the
proof of Claim 7. \qed
\enddemo

Now, let $v$ be the last letter of the word $V$. We consider two
cases.

\proclaim {Case I} $v^{\pm 1}$ occurs only once (at the end) in
the word $V$.
\endproclaim

Put $V \equiv V'tv$. If $v^{\pm 1} \notin \varPhi(U)$, then the
Whitehead automorphism $$\z=\z(\{v^{-1}, t\}, \, t)$$ fixes $x_1$
and $U$, but $\|\z(V)\|=|V'v|=|V|-1$. This contradicts the choice
of $V$. So $v^{\pm 1}$ must occur in $\varPhi(U)$. Besides, it
follows from $\overline{\b(V)} \equiv V_6e$ that $v \in {\Cal Q}$;
hence $v \in {\Cal Q} \cap \varPhi(U)$. Since ${\Cal T}_2 =
\emptyset$, $v \in {\Cal T}_1$. This implies that $y^{-1} \in
C_{x_1}(v, U)$, contrary to the existence of a Whitehead
automorphism $\a$ having the properties in (E1). Case I is
complete.

\proclaim {Case II} $v^{\pm 1}$ occurs more than once in the word $V$.
\endproclaim

Let us define
$$
\split
{\Cal L}_1=\{d \in {\Cal Q} \cap \varPhi(V)\,|\,&\text{there\ exists\ a\ path}\
q_d=d-v\ \text{of\ length}\ \ge 0\ \text{in}\ \varPhi(V)\\
&\text{that\ connects}\ d\ \text{with}\ v\ \text{and\ does\ not\
pass\ through}\ x_1^{\pm
1}\};\\
{\Cal L}_2=[{\Cal Q} \cap \varPhi(V)] \,\setminus \, {\Cal L}_1.
\endsplit
$$
If $e \in {\Cal L}_1$, then clearly $e \in C_{x_1}(v,V)$. Also, it
follows from $e \in {\Cal T}_1$ that $y^{-1} \in C_{x_1}(e, U)$.
Combining these facts gives us a contradiction to the existence
of a Whitehead automorphism $\a$ with the properties in (E1);
this contradiction enables us to have $$e \notin {\Cal L}_1.$$

Recall that $\overline{\b(V)} \equiv V_6e$. We treat two subcases
separately.

\proclaim {Case II.1} $V_6\not = V$.
\endproclaim

Consider the Whitehead automorphisms $$\rho_1=\rho_1({\Cal L}_1
\cup \{e\}, \, e) \quad \text {and} \quad \rho_2=\rho_2({\Cal L}_2
\cup \{e\}, \, e).$$ Since $e \notin {\Cal L}_1$ and $V_6\not = V$,
we see that $\|\rho_1(V)\| \ge |V|+2$. It then follows from the
observation
$$|V|+1 \ge \|\b(V)\|=\|\rho_1(V)\|+\|\rho_2(V)\|-|V|$$
that
$\|\rho_2(V)\| < |V|$. Here, assume that ${\Cal L}_1 \cap
\varPhi(U) = \emptyset$, and consider the Whitehead automorphism
$$\s=\s({\Cal T}_1 \cup {\Cal L}_2, \, e).$$ Since ${\Cal L}_1 \cap
\varPhi(U) = \emptyset$, ${\Cal L}_1 \cap ({\Cal T}_1 \cup {\Cal
L}_2) = \emptyset$. Also since ${\Cal T}_2 = \emptyset$, ${\Cal
Q} \cap \varPhi(U) = {\Cal T}_1 \subseteq {\Cal T}_1 \cup {\Cal
L}_2$. Thus we have $\s(U)=\b(U)$ and $\s(V)=\rho_2(V)$, so
that $\|\s(U)\|=|U|$ and $\|\s(V)\| < |V|$, contrary to the
choice of $V$ (note $\s(x_1)=x_1$ because $x_1^{\pm 1} \notin
{\Cal T}_1 \cup {\Cal L}_2 \subseteq {\Cal Q}$). So we must have
${\Cal L}_1 \cap \varPhi(U) \not = \emptyset$, say, $h \in {\Cal
L}_1 \cap \varPhi(U)={\Cal L}_1 \cap {\Cal T}_1$. But then we
have $h \in C_{x_1}(v,V)$ and $y^{-1} \in C_{x_1}(h, U)$,
contradicting the existence of a Whitehead automorphism $\a$
having the properties in (E1). Case II.1 is complete.

\proclaim {Case II.2} $V_6 = V$.
\endproclaim

This can happen only when
$${\Cal Q} \cap \varPhi(V) = \bigcup_{k \in
{\Cal Q} \cap \varPhi(V)}C_e(k, V). \tag 14
$$
It then follows that
$C_e(v, V) \subseteq {\Cal Q} \cap \varPhi(V)$, because $v \in
{\Cal Q} \cap \varPhi(V)$.

\proclaim {Case II.2.1} $C_e(v, V) \cap \varPhi(U) \not =
\emptyset$ .
\endproclaim

In this case, we have from ${\Cal Q} \cap \varPhi(U) =
{\Cal T}_1$ that $C_e(v, V) \cap {\Cal T}_1 \not = \emptyset$. If
$C_e(v, V) = C(v, V)$, then since $x_1^{\pm 1} \notin C_e(v, V)
\subseteq {\Cal Q}$, we have $C_e(v, V)=C(v, V)=C_{x_1}(v,
V)$, so that $C_{x_1}(v, V) \cap {\Cal T}_1 \not = \emptyset$.
But this yields a contradiction to the existence of a Whitehead
automorphism $\a$ with the properties in (E1). Thus,
$$C_e(v, V) \not = C(v, V).$$ It then follows that $e$ or $e^{-1}
\in C(v, V)$. If $e \in C(v, V)$, then there exists a path $q_e=e-v$
of length $\ge 0$ in $\varPhi(V)$ which connects $e$ with $v$.
Since $x_1^{\pm 1} \notin C_e(v, V) \subseteq {\Cal Q}$, the path
$q_e$ does not pass through $x_1^{\pm 1}$. But then, by
the definition of ${\Cal L}_1$, $e \in {\Cal L}_1$, contrary to $e \notin {\Cal L}_1$. Hence, we must have
$$e \notin C(v, V) \quad \text {and} \quad e^{-1} \in C(v,
V). \tag 15
$$

Define $${\Cal K}_1 = C_e(v, V) \cap {\Cal L}_1 \quad \text{and}
\quad {\Cal K}_2=C_e(v, V) \cap {\Cal L}_2.$$ In order to avoid a
contradiction to the properties of $\a$ in (E1), we must have ${\Cal K}_1 \cap {\Cal T}_1 \subseteq {\Cal L}_1 \cap {\Cal T}_1 =
\emptyset$. Here, consider the Whitehead automorphisms
$$\mu_1=\mu_1({\Cal K}_1 \cup \{e\}, \,e), \quad \mu_2=\mu_2({\Cal K}_2
\cup \{e\}, \,e) \quad \text{and} \quad \mu_3=\mu_3\bigl(({\Cal
L}_1 \cup {\Cal L}_2 \cup {\Cal T}_1) \setminus {\Cal K}_1, \, e
\bigl).$$ In view of $(15)$ and the definition of ${\Cal K}_1$,
we see that $\|\mu_1(V)\| \ge |V|+2$; hence it follows from the
observation
$$ |V|+1=|Ve|=\|\mu_1(V)\|+\|\mu_2(V)\|-|V|
\tag 16 $$ that $\|\mu_2(V)\| < |V|$. Since ${\Cal K}_1 \cap
{\Cal T}_1 = \emptyset$, we have ${\Cal Q} \cap
\varPhi(U)={\Cal T}_1 \subseteq [({\Cal L}_1 \cup {\Cal L}_2 \cup
{\Cal T}_1) \setminus {\Cal K}_1] \subseteq {\Cal Q}$; thus $\mu_3(U)=\b(U)$. Also,
in view of $(14)$ and the definitions of ${\Cal K}_1$ and ${\Cal
K}_2$, we can observe that $\mu_3(V)=\mu_2(V)$, contrary to the
choice of $V$ (note $\mu_3(x_1)=x_1$ because $x_1^{\pm 1} \notin
[({\Cal L}_1 \cup {\Cal L}_2 \cup {\Cal T}_1) \setminus {\Cal
K}_1] \subseteq {\Cal Q}$). Case II.2.1 is complete.

\proclaim {Case II.2.2} $C_e(v, V) \cap \varPhi(U) = \emptyset
\quad \text{and} \quad C_e(v, V)= C(v, V)$.
\endproclaim

Let $u$ be the first letter of the word $V$. Then
$$u^{-1} \notin C(v, V)=C_e(v, V);$$ for if $u^{-1} \in C(v, V)$ then it
would follow from $x_1^{\pm 1} \notin C(v, V)= C_e(v, V)
\subseteq {\Cal Q}$ that $u^{-1} \in C(v, V) = C_{x_1}(v, V)$,
contradicting $\a(V)=Vx_1^{-1}$ in (E1).

If $C(v,V)$ has a vertex of degree 1 that is different from $v$,
say, $w$, then the Whitehead automorphism $$\nu_1=\nu_1(\{w, c\},
\, c),$$ where $c$ is the adjacent vertex of $w$ in $\varPhi(V)$ (here $c$ cannot be $w^{-1}$), fixes $x_1$ and $U$ (because $w
\notin \varPhi(U)$ by the hypothesis of this case). But
$\|\nu_1(V)\| \le |V|-1 < |V|$. This contradicts the choice of
$V$, so every vertex in $C(v, V)$ different from $v$ must have
degree at least 2. This means that if $k \in C(v, V)$ then the
sum of occurrences of $k$ and $k^{-1}$ in $V$ must be at least
two.

Now assume that there is a vertex, say, $d$, in $C(v,V)$ such that
$d \in C(v,V)$ but ${d}^{-1} \notin C(v,V)$. Then since $x_1^{\pm
1} \notin C(v, V)= C_e(v, V) \subseteq {\Cal Q}$, for the
Whitehead automorphism
$$\nu_2=\nu_2\bigl(C(v, V),\, d),$$ we have $\nu_2(x_1)=x_1$. Also
since $C(v, V) \cap
 \varPhi(U)= C_e(v, V) \cap \varPhi(U) = \emptyset$ by the hypothesis of this case,
 $\nu_2(U)=U$. Furthermore, we see that if $d
\not = v$ then $\|\nu_2(V)\|=\|\pi_{d}(V)d\|$ (recall that $u^{-1}
 \notin C(v, V)$); whereas if $d=v$ then
$\|\nu_2(V)\|=\|\pi_{d}(V)\|$. Since ${d}^{\pm 1}$ occurs at
least twice in $V$, $\|\pi_{d}(V)\| \le |V|-2$; hence
$\|\nu_2(V)\| <|V|$. This contradicts the choice of $V$. Thus if
$d \in C(v,V)$, then ${d}^{-1}$ also has to occur in $C(v,V)$.
This yields that $C(v,V)=\varPhi(V)$, contradicting $u^{-1} \notin C(v, V)$. This completes Case II.2.2.

\proclaim {Case II.2.3} $C_e(v, V) \cap \varPhi(U) =
\emptyset \quad \text{and} \quad C_e(v, V) \not = C(v, V)$.
\endproclaim

Since $C_e(v, V) \not = C(v, V)$, we have, for the same reason as
in obtaining $(15)$, that $$e \notin C(v, V) \quad \text {and}
\quad e^{-1} \in C(v, V). \tag 17 $$ Define
$${\Cal M}_1 =
C_e(v, V) \cap {\Cal L}_1 \quad \text{and} \quad {\Cal
M}_2=C_e(v, V) \cap {\Cal L}_2,$$ and consider the Whitehead
automorphisms $$\xi_1=\xi_1({\Cal M}_1 \cup \{e\}, \, e) \quad
\text{and} \quad \xi_2=\xi_2({\Cal M}_2 \cup \{e\}, \, e).$$ By
$(17)$ and the definition of ${\Cal M}_1$, we see that
$\|\xi_1(V)\| \ge |V|+2$. Then it follows from an observation
similar to $(16)$ that $\|\xi_2(V)\|<|V|$. Also, since $x_1^{\pm
1} \notin {\Cal M}_2 \subseteq C_e(v, V) \subseteq {\Cal Q}$ and
${\Cal M}_2 \cap \varPhi(U) \subseteq C_e(v, V) \cap \varPhi(U) =
\emptyset$ by the hypothesis of this case, we have $\xi_2(x_1)=x_1$ and $\xi_2(U)=U$. A contradiction to the choice
of $V$ completes Case II.2.3.

The proof of the Theorem is now completed. \qed

\heading Acknowledgements 
\endheading

The author is greatly indebted to Professor S. V. Ivanov at the
University of Illinois at Urbana-Champaign for his valuable advice
and kind encouragement. She also thanks the referee for many helpful comments and suggestions. 

\heading References
\endheading

\roster

\item"[1]" G. Baumslag, A. G. Myasnikov and V. Shpilrain, {\it Open
problems in combinatorial group theory}, Contemp. Math. {\bf 250}
(1999), 1--27.

\item"[2]" S. V. Ivanov, {\it On endomorphisms of free groups that
preserve primitivity}, Arch. Math. {\bf 72} (1999), 92--100.

\item"[3]" R. C. Lyndon and P. E. Schupp, {\it Combinatorial Group Theory}, Springer-Verlag, New York/Berlin, 1977.

\item"[4]" J. McCool, {\it A presentation for the automorphism group
of a free group of finite rank}, J. London Math. Soc. {\bf 8}
(1974), 259--266.

\item"[5]" V. Shpilrain, {\it Recognizing automorphisms of the free groups},
Arch. Math. {\bf 62} (1994), 385--392.

\item"[6]" V. Shpilrain, {\it Generalized primitive elements of a
free group}, Arch. Math. {\bf 71} (1998), 270--278.

\item "[7]" J. H. C. Whitehead, {\it Equivalent sets of elements in a
free group}, Ann. of Math. {\bf 37} (1936), 782--800.

\endroster

\enddocument